\documentclass[12pt,spanish]{article}

\usepackage[utf8]{inputenc}
\usepackage{enumerate}
\usepackage{graphicx,graphics}
\usepackage[usenames]{color}
\usepackage{amsfonts}
\usepackage{amssymb}
\usepackage{amsthm}
\usepackage{dsfont}
\usepackage{eufrak}
\usepackage{MnSymbol,wasysym}
\usepackage{amsmath}
\usepackage[final]{pdfpages}
\newtheorem{theorem}{Theorem}
\newtheorem{definition}{Definition}
\newtheorem{lemma}{Lemma}
\newtheorem{proposition}{Proposition}
\newtheorem{corollary}{Corollary}
\newtheorem{Problem}{Problem}
\begin{document}

\title{Continuous selections, prime numbers and a covering type property}

\maketitle

\author{Jorge Antonio Cruz Chapital}

\begin{abstract}
Let $(X,\tau)$ be a Hausdorff space and $n\in\omega$. We prove that if $X$ admits a continuous selection over $\mathcal{F}_{n}(X)$ (nonempty subsets of $X$ of cardinality at most $n$), then for every $n\leq m\leq 2n$ such that $m$ is not a prime number, $X$ admits a continuous selection over $[X]^m$ (subsets of $X$ of cardinality $m$). As a consequence of this, a space $X$ admits a continuous selection for every natural number if and only if the same is true for every prime number. For Hausdorff spaces $(X,\tau)$ which admit continuous selections over $[X]^2$, we characterize the existence of continuous selections over $[X]^n$ for $n\geq 2$, in terms of a covering-type property.
\end{abstract}



\section{Preliminaries}
\label{sec1}
Our notation is fairly standard. $[X]^n$ is the set of all subsets of $X$ of cardinality $n$. If $(X,\tau)$ is a topological space, $\mathcal{F}(X)$ will be the set of all nonempty closed subsets of $X$, and if $n\in \omega$, then $F_{n}(X):=\{A\in \mathcal{F}(X)\; | \; |A|\leq n\}$. We can endow $\mathcal{F}(X)$ with the topology generated by  the sets$$\langle \mathcal{V}\rangle:=\{ A\in \mathcal{F}(X)\; | \; \forall V\in\mathcal{V}(A\cap V\not=\emptyset) \wedge A\subseteq \bigcup \mathcal{V}\}$$
where $\mathcal{V}$ is a finite subset of $\tau$. We will refer to this topology as the Vietoris topology $\tau_V$. All the spaces considered here are at least $T_2$, and from now on, we fix a space $(X,\tau)$.\\\\
The study of continuous selections began in the early 1950's by E. Michael \cite{Michael}. Since then, there has been an impressive amount of research that seems far from ending in the near future.  

\begin{definition}Given $\mathfrak{F}\subseteq \mathcal{F}(X)$, we say that $f:\mathfrak{F}\longrightarrow X$ is a selection if for every $A\in \mathfrak{F}$ we have that $f(A)\in A$. Moreover, $f$ is a continuous selection if it is a selection and it is continuous with respect to the Vietoris topology restricted to $\mathfrak{F}$. Additionally, given $n\in \omega$ we define \begin{itemize}
\item$Sel_n(X):=\{f:[X]^n\longrightarrow X\; | \; \textit{f is a selection}\}.$
\item$Sel_{\leq n}(X):=\{f:F_n(X)\longrightarrow X\; | \; \textit{f is a selection}\}.$
\item $Sel^c_n(X):=\{f:[X]^n\longrightarrow X\; | \; \textit{f is a continuous selection}\}.$
\item $Sel^c_{\leq n}(X):=\{f:F_n(X)\longrightarrow X\; | \; \textit{f is a continuous selection}\}$.\\
\end{itemize}
\end{definition}

Selecting and ordering are two closely related concepts. Examples of this can be found in \cite{Michael} where it is proven that connected spaces for which $Sel^c_2(X)\not=\emptyset$ are exactly those spaces which are weakly orderable, or in \cite{Mill}, where it is proven that the same is true for compact spaces. In \cite{Mill} it was asked if the same was true for locally compact spaces. A negative answer was given to that question in \cite{Hrusak} by constructing an almost disjoint family $\mathcal{A}$ for which $Sel^c_2(\Psi(\mathcal{A}))\not=\emptyset$ but $\Psi(\mathcal{A})$ is not weakly orderable. From this, we can conclude that although these concepts behave similar in many ways, they are not equivalent. Because of this, it is natural to ask how close are spaces for which   $Sel^c_2(X)\not=\emptyset$ to being weakly orderable. For example, if $(Y,\tau_y)$ is a weakly orderable space, then for every $n\in \omega$ we will have that $Sel^c_{\leq n}(Y)\not=\emptyset$. So in general, what can we say of $Sel^c_n(X)$, knowing that $Sel^c_2(X)\not=\emptyset$?. The following question was posed in \cite{Guno}.

\begin{Problem}[\cite{Guno}]\label{Problem1}Does there exist a space $X$ such that $Sel_2^c(X)\not=\emptyset$, but $Sel^c_{\leq n}(X)=\emptyset$ for some $n>2$?

\end{Problem}
This question is still open even in the case $n=3$. In \cite{Garcia}, it is proven that if $Sel^c_2(X)\not=\emptyset$, then it is equivalent that $Sel^c_3(X)\not=\emptyset$ and $Sel^c_{\leq 3}(X)\not=\emptyset$. This result was later generalized in \cite{Gutev} by the following Theorem.

\begin{theorem}[\cite{Gutev}]\label{Gutev}Let $n\in \omega$ such that $Sel^c_{\leq n}(X)\not=\emptyset$. The following statements are equivalent:
\begin{enumerate}[a)]
\item $Sel^c_{\leq{n+1}}(X)\not=\emptyset$ .
\item $Sel^c_{n+1}(X)\not=\emptyset$.
\end{enumerate}
\end{theorem}
By means of Theorem \ref{Gutev}, we can restate  Problem \ref{Problem1} as:

\begin{Problem}\label{Problem1.5}Does there exist a space $X$ for which $Sel^c_2(X)\not=\emptyset$, but $Sel^c_n(X)=\emptyset$ for some $n>2$?
\end{Problem}

A partial answer to this question can be found \cite{Hrusak}, where the next Theorem shows up. 

\begin{theorem}[\cite{Hrusak}] Let $X$ be a separable space for which $Sel^c_2(X)\not=\emptyset$. Then, there exists an orderable space $L$ and a continuous map $f:X\longrightarrow L$ such that $|f^{-1}[\{y\}]|\leq 2$ for every $y\in L$. In particular, $X$ admits a continuous selection over $[X]^{<\omega}$.
\end{theorem}

In \cite{Guno3} it is proven that $Sel^c_4(X)\not=\emptyset$ provided that $Sel^c_{\leq 3}(X)\not=\emptyset$. This result was later generalized in \cite{Gutev} by means of the following Theorem.
\begin{theorem}[\cite{Gutev}]\label{2n+1implica2n+2} Supose that $X$ is a  space for which $Sel^c_{2n+1}(X)\not=\emptyset$ for some $n\in \omega$. Then $Sel^c_{2n+2}(X)\not=\emptyset$.
\end{theorem}
So in particular, if Problem \ref{Problem1.5} has an affirmative answer, there must exist $X$ for which $Sel^c_2(X)\not=\emptyset$ but $Sel^c_n(X)\not=\emptyset$ for some odd number $n$. Moreover, if we aim to answer Problem \ref{Problem1.5}, we may only worry about $Sel^2_n(X)$ for odd numbers $n$.\\

The aim of this paper is to generalize Theorem \ref{2n+1implica2n+2} by means of Theorem $\ref{teoprin1}$, which will be achieved in Section \ref{sec1}. In Section \ref{sec2} we will define the concept of good chains, and use it characterize spaces $X$ for which $Sel^c_n(X)\not=\emptyset$ provided that $Sel^c_2(X)\not=\emptyset$.\\

Let us start with some notation.

\begin{definition}Let $\mathfrak{F}\subseteq F(X)$, $f:\mathfrak{F}\longrightarrow X$ a selection, and $\vec{x}\in Dom(f)$.  We will write $$\vec{x}\longrightarrow_f y$$ if $f(\vec{x})=y$. Additionally, if $\mathcal{U}\subseteq P(X)$, then We will write $$\mathcal{U}\rightrightarrows_f V$$ if $V\in \mathcal{U}$ and for every $\vec{x}\in \langle \mathcal{U}\rangle\cap Dom(s)$, it happens that $s(\vec{x})\in V$. If $f$ is clear from context then we'll omit it's reference in the notation.
\end{definition}
By means of this Definition, we intend to think about continuous selections as directed graphs.
The next Proposition is easy to prove, and characterizes the continuity of selections in a way coherent with our intuition.
\begin{proposition}\label{folklore}Let $\mathfrak{F}\subseteq F(X)$, and $s:\mathfrak{F}\longrightarrow X$ a selection. The following statements are equivalent:
\begin{enumerate}[a)]
\item $f$ is continuous.
\item for every $n\in \omega$ and $\{x_1,\dots,x_n\}\in Dom(s)\cap[X]^n$, there exist $U_1,\dots,U_n\in \tau$ pairwise disjoint open neighborhoods of $x_1,\dots,x_n$  respectively, , y $j\leq n$, such that $\{U_1,\dots,U_n\}\rightrightarrows U_j$.\\
\end{enumerate}
\end{proposition}

\begin{definition}Let $n\leq m\in \omega$, $f\in Sel_n(X)$ and $\mathcal{U}\in [\tau\backslash\{\emptyset\}]^m$ a set of pairwise disjoint open sets.  We will say that  $\mathcal{U}$ preserves $f$-relations if for every $\mathcal{V}\in [\mathcal{U}]^n$ there exists some $V\in \mathcal{V}$ such that  $$\mathcal{V}\rightrightarrows V.$$ 
In general, if $k,m\in \omega$, $f\in Sel_{\leq k}(X)$ and $\mathcal{U}\in [\tau\backslash\{\emptyset\}]^m$ a set of pairwise disjoint open sets. We will say that $\mathcal{U}$ preserves $f$-relations if for every $0<n\leq min(k,m)$, it happens that $\mathcal{U}$ preserves $f|_{[X]^n}$-relations.
\end{definition}
Proposition \ref{folklore} told us, that by  continuity, we can fatten the arrows between points to make them arrows between open sets. Next proposition tells us that we can do the same if we consider "finite subgraphs" of $X$.
\begin{proposition}\label{presrel}Let $n\in \omega$ and $f\in Sel^c_{\leq n}(X)$(or $f\in Sel^c_n(X)$). For every $n\leq m\in \omega$ and for every $\{x_1,\dots x_m\}\in [X]^m$ There exists $U_1,\dots,U_m\in \tau$ pairwise disjoint open neighborhoods of $x_1,\dots,x_m$ respectively, such that $\{U_1,\dots,U_m\}$ preserves $f|_[X]^k$-relations for every $k\leq n$.
\end{proposition}

\section{Isomorphisms of Selections}
\label{sec2}
The main purpose of this section is to prove Theorem \ref{teoprin1}. Before explaining the main idea of the proof, we will give some definitions.
\begin{definition}Given $n\in \omega$, $s\in Sel_n(Z)$, y $r\in Sel_n(Y)$. We will say that $f$ and $g$ are isomorphic if there exists $\phi: Z\longrightarrow Y$ biyective such that for every $\vec{z}\in [Z]^n$ it happens that 
 $$\phi[\vec{z}]\longrightarrow_g \phi(f(\vec{z})).$$ In this case, we will say that $\phi$ is an isomorphism between f and g.
\end{definition}

\begin{definition}Let $n\leq k\in \omega$, $f\in Sel_{\leq k}(X)$ and $g\in Sel_n(Y)$. We define $$\mathcal{P}(g):=\{\vec{x}\in [X]^m\; | \; f|_{[\vec{x}]^n} \; \textit{ is isomorphic to } g\}.$$
\end{definition}

Now, let us explain the strategy. Suppose that we have a space $X$ for which $Sel^c_{\leq k}(X)\not=\emptyset$ for some $k\in \omega$, and let $f$ be a continuous selection that testifies this fact. For $m$ between $k$ and $2k$, and $n\leq k$ we will partition $[X]^m$ into $\{\mathcal{P}(g)\; | \; g\in Sel_n(m)\}$ and show that each piece of the partition is clopen in $[X]^m$. Having done this, the problem will be reduce to find a continuous selection over $\mathcal{P}(g)$ for every $g\in Sel_n(m)\}$. This will be done in two steps; First we will show that if $g$ is not too "messy", $\mathcal{P}(g)$ admits a continuous selections, and then we will show that if $n$ is prime and divides $m$, then there are no "messy" $g's$.\\

In the next two propositions we will show that each $\mathcal{P}(g)$ is clopen.

 \begin{proposition}\label{isoabiertos} Let $n\leq m\in\omega$, $f\in Sel_n(X)$ and $\mathcal{U}=\{U_1,\dots, U_{m}\}\in [\tau\backslash\{\emptyset\}]^m$ who's elements are pairwise disjoint. If $\mathcal{U}$ preserves relations, then for every $\vec{x}=\{x_1,\dots,x_{m}\},\vec{y}=\{y_1,\dots,y_m\}$ such that for every $i\leq m$ it happens that $x_i,y_i\in U_i$, the function  $\phi^{\vec{x}}_{\vec{y}}(f):\vec{x}\longrightarrow \vec{y}$ given by  $$\phi^{\vec{x}}_{\vec{y}}(s)(x_i)=y_i$$
 
is an isomorphism between $f|_{[\vec{x}]^n}$ a $f|_{[\vec{y}]^n}$.

\begin{proof}Let $i_1<\dots<i_{n}\leq m$ and $j\leq n$ such that $$\{x_{i_1},\dots,x_{i_{n}}\}\longrightarrow x_j.$$ Since $\mathcal{U}$ preserves relations, we can conclude that $\{U_{i_1},\dots,U_{i_{n}}\}\rightrightarrows U_j$, and since $y_{i_1}\in U_{i_1},\dots,y_{i_{n}}\in U_{i_{n}}$ then $$\phi^{\vec{x}}_{\vec{y}}(f)[\{x_{i_1},\dots,x_{i_{n}}\}]=\{y_{i_1},\dots, y_{i_{n}}\}\longrightarrow y_{i_j}=\phi^{\vec{x}}_{\vec{y}}(f)(x_{i_j}).$$

Consequently $\phi^{\vec{x}}_{\vec{y}}(f)$ is isomoprhism.
\end{proof}
\end{proposition}

 \begin{proposition}\label{clopenp(g)} Let $n\in \omega$, $Y$ such that $|Y|=m$ , $f\in Sel^c_n(X)$ y $g\in Sel_n(Y)$. Then $\mathcal{P}(g)$ is clopen in $[X]^m$.
 \begin{proof}Since $\{\mathcal{P}(g)\; |\; g\in Sel_n(Y)\}$ form a partition of $[X]^m$, it's enough to proof that for every $r\in Sel_n(Y)$, $\mathcal{P}(r)$ is open.\\\\
 Take $g\in Sel_n(Y)$ and $\vec{x}=\{x_1,\dots,x_m\}\in \mathcal{P}(g)$. Using Proposition \ref{presrel}, we find $U_1,\dots, U_m\in \tau$ pairwise disjoint open neighborhoods of $x_1,\dots, x_m$ respectively, such that $\{U_1,\dots,U_n\}$ preserves $f$-relations. To finish, just notice that $\vec{x}\in \langle U_1,\dots, U_m \rangle$ and by Proposition \ref{isoabiertos} it
 follows that $ \langle U_1,\dots, U_m\rangle\subseteq \mathcal{P}(g)$.
 \end{proof}
\end{proposition}
By now, our first task is completed. The next definition is crucial to understand which are the messy $g$'s that we were referring to.
\begin{definition} Given $n\in \omega$, $g\in Sel_n(Y)$ and $x\in Y$, we define $$\mathcal{W}_g(x):=|g^{-1}[\{x\}]|.$$ Additionally, if $k\in \omega$ we define $$\mathcal{Q}_g(k):=\{x\in Y\; | \; \mathcal{W}_g(x)=k\}.$$ If $g$ is clear by context, we'll omit it.
\end{definition}
\begin{lemma}\label{sel2p(g)} Let $k,n,m\in \omega$ with $n\leq k$, and $Y$ such that $|Y|=m$, $f\in Sel^c_{\leq_k}(X)$ y $g\in Sel_{n}(Y)$. If there exists $x,y\in Y$ such that $\mathcal{W}_g(x)\not=\mathcal{W}_g(y)$  and $\frac{m}{2}\leq k$ then there exists a continuous selection $h$ over $\mathcal{P}(r)$.
\begin{proof} The set $\mathcal{Q}=\{\mathcal{Q}_g(r)\; | r\in \omega \textit{ and } 
 \mathcal{Q}_g(r)\not=\emptyset\}$ is a partition of $Y$, and by hipotesis, it has at least two elements. Thus, there exists $r\in \omega $ such that $ \mathcal{Q}_g(r)\not=\emptyset $ y $|\mathcal{Q}_g(r)|\leq \frac{m}{2}$. take $r_0\in \omega$  the least natural with this property and let $k_0:=|\mathcal{Q}_g(r_0)|$.\\
 
We define $h:\mathcal{P}(g)\longrightarrow X$ given by $h(\vec{x})=f(Q_{f|_{[x]^n}}(r_0))$. We know $h$ is well defined, since for every $\vec{x}\in \mathcal{P}(g)$ it happens that $|\mathcal{Q}_{f|_{[\vec{x}]^n}}(r_0)|=k_0\leq k$. To proof that $h$ is continuous, take $\vec{x}=\{x_1,\dots,x_m\}\in \mathcal{P}(g)$. Applying \ref{presrel}, we can find $U_1,\dots, U_m\in \tau$ pairwise dijsoint open neighborhoods of  $x_1,\dots,x_m$ respectively, such that $\{U_1,\dots, U_m\}$ preserves $f|_{[X]^n}$-relations  and $f|_{[X]^{k_0}}$-relations. \\\\ Consider  $i_1<\dots <i_{k_0}\leq n$ such that $\mathcal{Q}_{f|_{[x]^n}}(r_0)=\{x_{i_1},\dots,x_{i_{k_0}}\}$, thus for every $\vec{y}=\{y_1,\dots,y_n\}$ with $y_1\in U_1,\dots, y_m\in U_m$ it happens that
 $$\{y_{i_1},\dots,y_{i_{k_0}}\}=\mathcal{Q}_{f|_{[\vec{y}]^n}}(r_0)=\phi^{\vec{x}}_{\vec{y}}[\mathcal{Q}_{f|_{[\vec{x}]^n}}(r_0)].$$
Since $\{U_1,\dots, U_m\}$ preserves $f|_{[X]^{k_0}}$-relations, so it does $\{U_{i_1},\dots, U_{i_{k_0}}\}$. Notice that for every $j\leq k_0$ we have that $x_{i_j},y_{i_j}\in U_{i_j}$, consequently $$h(x)=f(Q_{f|_{[\vec{x}]^n}}(r_0))\in U_{i_j} \textit{ if and only if } h(y)=f(Q_{f|_{[\vec{y}]^n}}(r_0))\in U_{i_j}.$$
We conclude that $\{U_1,\dots, U_m\}\rightrightarrows _h U_j$ where $j\leq n$ is such that $\vec{x}\longrightarrow x_j$, thus $h$ is continuous.
\end{proof}
\end{lemma}

 Basically, Lemma \ref{sel2p(g)} is telling us that the messy $g$'s are exactly those for which $\mathcal{W}_g$  is constant. We will study this in more detail on section \ref{secgoodchains}.

\begin{lemma}\label{extsel} let $k,m\in \omega$ and $f\in Sel^c_{\leq_k}(X)$. If $\frac{m}{2}\leq k$, and there is $n\leq k$ such that for every $Z\in [X]^m$ there exists $x,y\in Z$ such that $\mathcal{W}_{f|_{[Z]^n}}(x)\not=\mathcal{W}_{f|_{[Z]^n}}(y)$, then $Sel^c_m(X)\not=\emptyset$. 
\begin{proof} Fix a set $Y$ such that $|Y|=m$. Then for every $g\in Sel_n(Y)$ such that $\mathcal{P}(g)\not=\emptyset$, The hipotesis of the Lemma \ref{sel2p(g)} are fulfilled , thus there exists a continuous selection $h_g:\mathcal{P}(g)\longrightarrow X$. Since $\mathcal{P}(g)$ is clopen, then $$h=\bigcup\{h_g\; |\; g\in Sel_n(Y) \textit{ y } \mathcal{P}(g)\not=\emptyset\}$$
is a continuous selection over $[X]^m$.
\end{proof}
\end{lemma}

\begin{lemma}\label{prrimos}Let $p\leq m\in \omega$ and suppose that $p$ is prime, and $|Y|=m$. If there exists $g\in Sel_p(Y)$such that for every $x,y\in Y$ if happens that $\mathcal{W}(x)=\mathcal{W}(y)$, then $p$ does not divides $m$.
\begin{proof}Le $g$ be as in the hypothesis. Suppose toward a contradiction that there is $j\in \omega$ such that $pj=m$, and let $k$ such that for every  $x\in Y$ if happens that $k=\mathcal{W}(x)$. Notice that $$p(jk)=mk=\sum_{x\in Y}\mathcal{W}(x)=|[Y]^p|=\frac{m!}{p!(m-p)!}$$

Thus, $p$ divides $\frac{m!}{p!(m-p)!}$ ,but this is impossible since $p$ is prime and $p$ divides $m$. We conclude that $p$ does not divides $m$.
\end{proof}
\end{lemma}

\begin{theorem}\label{teoprin1} Let $p,k,m\in \omega$ and supose $ Sel^c_{\leq k}(X)\not=\emptyset$. If $p, \frac{m}{2}\leq k$ and $p$ divides $m$, then $Sel^c_m(X)\not=\emptyset$.

\begin{proof}Let $f\in Sel^c_{\leq k}(X)$Since $p$ divides $m$, applying Lemma \ref{prrimos} we can conclude that for every $Z\in [X]^p$  there exist $x,y\in Y$ such that $\mathcal{W}_{f|_{[Z]^n}}(x)\not=\mathcal{W}_{f|_{[Z]^n}}(y)$, thus, $k,m$ and $f$ fulfill the hypothesis of the Lemma  \ref{extsel}, and in consequence $Sel^c_m(X)\not=\emptyset$.
\end{proof}
\end{theorem}

\begin{corollary}\label{compuesto}Let $2\leq n$ and suppose that $ Sel^c_{\leq n}(X)\not=\emptyset$. If $n+1$ is not prime, then $Sel^c_{n+1}(X)\not=\emptyset$.
\end{corollary}

\begin{theorem}Let $(X,\tau)$ be a $T_2$ topological space. The following statements are equivalent:
\begin{enumerate}[a)]
\item For every $1<n\in\omega$,  $Sel^c_{n}(X)\not=\emptyset$.
\item For every prime $p$,  $Sel^c_p(X)\not=\emptyset$.
\end{enumerate}

\begin{proof}$a)\rightarrow b)$ trivial.
	\\\\$b)\rightarrow a)$ By induction. Let $1<n\in \omega$ and suppose that for every $1<m<n$ we know that $Sel^c_n(X)\not=\emptyset$. If $n$ is prime, we use the hypothesis . If $n$ is not prime, then we make use of Theorem \ref{Gutev} to conclude that $Sel^c_{\leq n-1}(X)\not=\emptyset$. Thus, by Corollary  \ref{compuesto} it follows that $Sel^c_n(X)\not=\emptyset$.\end{proof}
\end{theorem}

\section{Good chains}\label{secgoodchains}
The purpose of this section is to characterize the spaces $X$ for which $Sel^c_{n+1}(X)\not=\emptyset$ provided that $2\leq n$ and $Sel^c_n(X)\not=\emptyset$. To do this, recall that we showed in Lemma \ref{sel2p(g)} that for $k \leq n$ and $g\in Sel_k(n+1)$, we have that $\mathcal{P}(g)$ admits a continuous selection unless  $\mathcal{W}_g(x)=\mathcal{W}_g(y)$ for every $x,y\in n+1$. Now we'll fix our atention in the case $k=2$.

\begin{definition} Let $2\leq n\in \omega$ odd. We define $$Reg_n:=\{g\in Sel_2(n)\; |\; \forall x,y\in n\big(\mathcal{W}_g(x)=\mathcal{W}_g(y)\big)\}$$ 
   More over,  if $2\leq m\in \omega$ and $f\in Sel^c_{\leq m}(X)$, we define $$\mathcal{P}(Reg_n)=\bigcup_{g\in Reg_n}\mathcal{P}(g).$$
\end{definition}\label{regsel}

If $Sel^c_{\leq_n}(X)\not=\emptyset$ and $f$ testifies this fact, then $Sel^c_{n+1}(X)\not=\emptyset$ if and only if there is a continuous selection over $\mathcal{P}(Reg_{n+1})$. A simple argument shows that when $X$ is infinite then $[X]^{n+1}\backslash\mathcal{P}(Reg_{n+1})\not=\emptyset$, so analysing  $\mathcal{P}(Reg_{n+1})$ really simplifies the problem. \\

An important property of elements in  $Reg_n$ is the following.

\begin{proposition} Let $2\leq n\in \omega$ odd and $g\in Reg_n$. For every $x,y\in n$ such that $x\longrightarrow y$, there is $z\in n\backslash \{x,y\}$ such that $y \longrightarrow z\longrightarrow x.$
\begin{proof}Let $A_x=\{z\in n\backslash{y}\; | \; z\longrightarrow x\}$ and $A_y=\{z\in n\backslash{x}\; | \; z\longrightarrow y\}$. Since $g\in Reg_n$,  $|A_x|=\frac{n-1}{2}$ y $|A_y|=\frac{n-1}{2}-1$. Thus, $A_x\backslash A_y\not=\emptyset$. Notice that every $z\in A_x\backslash A_y$ works.\\
\end{proof}
\end{proposition}
Before explaining the main idea behind this characterization, let's write one more definition.
\begin{definition}For $m\in\omega$ we can define $$\mathfrak{F}_m(X):=\{\mathcal{U}\in [\tau\backslash \{\emptyset\}]^m\; |\; \forall U,V\in \mathcal{U}, \textit{ if } U\not=V \textit{ then } U\cap V=\emptyset \}.$$
\end{definition}
 
Suppose that $f\in Sel_{\leq n}(X)$ and we want to construct a continuous selection over $\mathcal{P}(Reg_{n+1})$, namely $g$ . For every $\vec{x}\in \mathcal{P}(Reg_{n+1})$ we may take $\mathcal{U}_x\in \mathfrak{F}_m(X)$ such that $\vec{x}\in \langle \mathcal{U}_x\rangle$ and define the value of the selection $g$ in such way that $\mathcal{U}\longrightarrow_g U$ for some $U\in \mathcal{U}_x$, in order to ensure continuity by Proposition \ref{folklore}. The problem is that in some cases, the "information" that caries some $\mathcal{U}$ about the selection $g$ is transferred to some other $\mathcal{V}$'s, in the sense that if we define the value of $g$ at $\langle \mathcal{U}\rangle$ then the value of those $\mathcal{V}'s$ will be determined as well, and in some cases, this fact may translate into contradictions if we define our $g$ in careless way.\\

The next few definitions capture the core of this problem. 

\begin{definition}\label{interbonitadef} Let $n\in \omega$ and  $\mathcal{U},\mathcal{V}\in \mathfrak{F}_n(X)$. We'll write $\mathcal{U}\smile \mathcal{V}$ if for every $U\in \mathcal{U}$ there exists a unique $V\in \mathcal{V}$ such that $U\cap V\not=\emptyset$. If this is the case, we'll denote as $\gamma_{(\mathcal{U},\mathcal{V})}:\mathcal{U}\longrightarrow \mathcal{V}$ the only function which satisfies that for every  $U\in \mathcal{V}$,  $\gamma_{(\mathcal{U},\mathcal{V})}(U)\cap U\not=\emptyset$.
 \end{definition}
 
Definition \ref{interbonitadef} is telling us the easiest way possible in which the value of a selection at $\mathcal{U}$  can determine the value of the selection at $\mathcal{V}$(under circumstances that wil be clear in the proof of Theorem \ref{teoprin2}). Now we define the concept of chain, which will describe essentially all the ways in which the information of a given selection can be passed from one set to another.
 
 \begin{definition}Let $n,m\in \omega$, $\mathfrak{D}\subseteq \mathfrak{F}_n(X)$ and $\mathcal{U}_0,\dots,\mathcal{U}_m\in \mathfrak{F}_n(X)$. We say that $\mathbb{K}=(\mathcal{U}_0,\dots,\mathcal{U}_m)$ is a $\mathfrak{D}$-chain from $\mathcal{U}_0$ to $\mathcal{U}_m$(or simply a $\mathfrak{D}$-chain), if for every $i<m$ it happens that $\mathcal{U}_i\in \mathfrak{D}$ and $\mathcal{U}_i\smile \mathcal{U}_{i+1}$. Additionally, we define the function  $$\Gamma_\mathbb{K}:=\gamma_{(U_{m-1},U_m)}\circ \dots \circ \gamma_{(U_0,U_1)}.$$

Finally, we'll say that $\mathcal{U}, \mathcal{V}\in  \mathfrak{D}$ are $\mathfrak{D}$-chainable if there exists a $\mathfrak{D}$-chain from  $\mathcal{U}$ to $\mathcal{V}$.\\
\end{definition}

\begin{definition}Let $n\in \omega$ and $\mathfrak{D}\subseteq \mathfrak{F}_n(X)$.  We say that $\mathfrak{D}$ is nice if for every $\mathcal{A},\mathcal{B}\in \mathfrak{D}$ the following statements are true:\begin{itemize}
	\item If $\langle \mathcal{A}\rangle \cap \langle{B}\rangle\not=\emptyset$ then $\mathcal{A}\smile \mathcal{B}$.
	\item If $\mathcal{A}$ and $\mathcal{B}$ are $\mathfrak{D}$-chainable, then for every $\mathbb{K}, \mathbb{L}$ \break $\mathfrak{D}$-chains from $\mathcal{A}$ to $\mathcal{B}$, we have that $\Gamma_\mathbb{K}=\Gamma_\mathbb{L}$.
	
\end{itemize}
	
\end{definition}

Before we state our characterization, let's give one last Lemma.
\begin{lemma}\label{interbonita} Let $2\leq n\in \omega$ even, $f\in Sel^c_{2}(X)$ and $\mathcal{U}$, $\mathcal{V}\in \mathfrak{F}_{n+1}(X)$ which preserve relations. If $\langle \mathcal{U}\rangle,\langle \mathcal{V}\rangle\subseteq \mathcal{P}(Reg_{n+1})$, then $\langle \mathcal{U}\rangle\cap\langle \mathcal{V}\rangle\not=\emptyset$ implies that $\mathcal{U}\smile \mathcal{V}$. \begin{proof} 
Suppose towards a contradiction, that there is  $U\in \mathcal{U}$ such that for every $V\in \mathcal{V}$ it happens that $U\cap V=\emptyset$. Notice that for every $x\in\langle\mathcal{V}\rangle$ we have that $x\cap U=\emptyset$. Thus $x\notin \langle\mathcal{U} \rangle$, and then $\langle \mathcal{U}\rangle\cap\langle \mathcal{V}\rangle=\emptyset$ which is impossible . \\\\To prove uniqueness, suppose towards a contradiction that there is $U\in \mathcal{U}$ and $V_0, V_1\in\mathcal{V}$ distinct such that for every $i\in 2$ it happens that $U\cap V_i\not=\emptyset$. Let $V_2,\dots, V_m\in \tau_\backslash\{\emptyset\}$ such that $\mathcal{V}=\{V_0,V_1,V_2,\dots, V_m \}$, y suppose without loss that $V_0\rightrightarrows V_1$. Since $\mathcal{V}$ preserves relations and  $\langle \mathcal{V}\rangle \subseteq\mathcal{P}(Reg_{n+1})$, it follows that there $g\in Reg_{n+1}$ such that $\langle\mathcal{V}\rangle\subseteq \mathcal{P}(g)$. Thus, making us of Proposition \ref{regsel} we find $2\leq j<n+1$ such that $V_1\rightrightarrows V_j\rightrightarrows V_0$. We divie the rest of the proof in two cases:\\

{\bf Case 1)} (There exists $U\not=U'\in \mathcal{U}$ such that $U'\cap V_j\not=\emptyset$) In this case, take $x_0\in U\cap V_0$, $x_1\in U\cap V_1$ y $x_2\in U'\cap V_j$. In one hand, since $V_1\rightrightarrows V_j\rightrightarrows V_0$ then $$x_1\longrightarrow x_2\longrightarrow x_0,$$

By the other side, since $\mathcal{U}$ preserves relations, $U\not=U'$, $x_2\in U'$ y $x_0,x_1\in U$, we have that $x_1\longrightarrow x_2$ if and only if $x_0\longrightarrow x_2$, which is a contradiction.\\

{\bf Case 2)} (For every $U\not=U'\in \mathcal{A}$ it happens that $U'\cap V_j=\emptyset$) Notice that in this case $V_j\cap U\not=\emptyset$. If for  every $U\not=U'\in \mathcal{A}$ we had that $U'\cap V_0=U'\cap V_1=\emptyset$ We could argue in a similar way as in the first paragraph of this proof, to conclude that $\langle \mathcal{U}\rangle\cap\langle \mathcal{V}\rangle=\emptyset$, Which is impossible. Thus, there is $i<2$ and $U\not=U'\in \mathcal{U}$ such that $U'\cap V_i\not=\emptyset$. Let $x_0\in U\cap V_{1-i}$, $x_1\in U\cap V_j$ and $x_2\in U'\cap V_i$. Proceeding in a completely  similar way to case 1 we get to the desired contradiction.
\end{proof}
\end{lemma}

Now, we are ready.
\begin{theorem}\label{teoprin2} Let $(X,\tau)$ be a topological space, $2\leq n\in \omega $ even and $f\in Sel^c_{\leq n}(X)$. The following statements are equivalent:

\begin{enumerate}[a)]
\item $Sel^c_{n+1}(X)\not=\emptyset$	.

\item There is a nice $\mathfrak{D} \subseteq \mathfrak{F}_{n+1}(X)$, such that $$\bigcup_{\mathcal{U}\in \mathcal{D}}\langle \mathcal{U}\rangle=\mathcal{P}(Reg_{n+1}).$$
\end{enumerate}
\begin{proof}
{$\bf{ \big(a)\Rightarrow b)\big)}$} Let $h\in Sel^c_{n+1}(X)$ and $g=f\cup h$. For every $\vec{x}=\{x_0,\dots,x_n\}\in \mathcal{P}(Reg_{n+1})$, let's take $U^{\vec{x}}_0,\dots, U^{\vec{x}}_n\in \tau\backslash \{\emptyset\}$ open pairwise disjoint neighborhoods of $x_0,\dots,x_n$ respectively,such that for every $i\leq n+1$ it happens that$\mathcal{U}^{\vec{x}}=\{U^{\vec{x}}_0,\dots, U^{\vec{x}}_n\}$ preserves $g|_{[X]^i}$-relations ( We can do this by applying repeatedly Lemma \ref{presrel} and intersecting witnesses ). Let $$\mathfrak{D}=\{\mathcal{U}^{\vec{x}}\;| \; \vec{x}\in \mathcal{P}(Reg_{n+1})\}. $$

Evidently $\mathfrak{D}\subseteq \mathfrak{F}_{n+1}(X)$. Since for every $\vec{x}\in \mathcal{P}(Reg_{n+1})$ it happens that $\mathcal{U}^{\vec{x}}$ preserves $g|_{[X]^2}$-relations and $ \vec{x}\in\langle \mathcal{U}^{\vec{x}}\rangle$, then $$\bigcup_{\mathcal{U}\in \mathcal{D}}\langle \mathcal{U}\rangle=\mathcal{P}(Reg_{n+1}).$$

 It's only left to show that $\mathfrak{D}$ is nice. For this, we take $\mathcal{U},\mathcal{V}\in \mathfrak{D}$, and prove the following:
 \begin{enumerate}
 \item	{ ( If $\langle \mathcal{U}\rangle \cap \langle \mathcal{V}\rangle\not=\emptyset$ then $\mathcal{U}\smile \mathcal{V}$.)} This is a direct consequence of Lemma \ref{interbonita}.
 \item {( If $\mathcal{U}$ and $\mathcal{V}$ are $\mathfrak{D}$-chainable, then for every $\mathbb{K}, \mathbb{L}$ $\mathfrak{D}$-chains from $\mathcal{U}$ to $\mathcal{V}$,it happens that $\Gamma_\mathbb{K}=\Gamma_\mathbb{L}$.)} To prove this, we'll make use of the next claim.\\\\
 {\bf \underline{Claim 1:}} if $\vec{x}=\{x_0,\dots,x_n\},\vec{y}=\{y_0,\dots,y_n\}\in \mathcal{P}(Reg_{n+1})$ are such that $\mathcal{U}^{\vec{x}}\smile \mathcal{U}^{\vec{y}}$, then for every $k\leq n$ and $i_0,\dots,i_k\leq n$ it happens that $$\{U^{\vec{x}}_{i_0},\dots,U_{i_k}\}\rightrightarrows U^{\vec{x}}_{i_0}\textit{ if and only if } \{\Gamma_{\mathbb{K}}(U^{\vec{x}}_{i_0}),\dots,\Gamma_{\mathbb{K}}(U^{\vec{x}}_{i_k})\}\rightrightarrows \Gamma_{\mathbb{K}}(U^{\vec{x}}_{i_0})$$
   where $\mathbb{K}=(U^{\vec{x}},U^{\vec{y}})$.\\
   
To prove this claim, notice that if $k\leq n$ and $i_0,\dots, i_k\leq n$, then for every $s\leq k$ we can choose $z_s\in U_{i_s}^{\vec{x}}\cap \Gamma_{\mathbb{K}}(U_{i_s}^{\vec{x}})$, and since  $\mathcal{U}^{\vec{x}}$ and $ \mathcal{U}^{\vec{y}}$ preserve $g|_{[X]^k}$- relations, then \break $\{U^{\vec{x}}_{i_0},\dots,U_{i_k}\}\rightrightarrows U^{\vec{x}}_{i_0}$ if and only if $ \{z_0,\dots, z_k\}\longrightarrow z_0$ if and only if \break $\{\Gamma_{\mathbb{K}}(U^{\vec{x}}_{i_0}),\dots,\Gamma_{\mathbb{K}}(U^{\vec{x}}_{i_k})\}\rightrightarrows \Gamma_{\mathbb{K}}(U^{\vec{x}}_{i_0}).$\\

Having proved Claim 1, let's continue with the proof of (2). Let $\vec{x}=\{x_0,\dots,x_n\},\vec{y}=\{y_0,\dots,y_n\}\in \mathcal{P}(Reg_{n+1})$ and suppose that $\mathcal{U}^{\vec{x}}$ and $\mathcal{U}^{\vec{y}}$ are chainable. Without loss, we can suppose that $x_0,\dots,x_n$ and $y_0,\dots, y_n$ are enumerated in such way that for every $i\leq n$ it happens that $\{x_i,\dots, x_n\}\longrightarrow x_i$ and $\{y_i,\dots,y_n\}\longrightarrow y_i$. Let $\mathbb{K}$ be a $\mathfrak{D}$-chain from $\mathcal{U}^{\vec{x}}$ to $\mathcal{U}^{\vec{y}}$ and let $k\in \omega$ and $\vec{z}_0,\dots,\vec{z}_k
\in \mathcal{P}(Reg_{n+1})$ be such that $\mathbb{K}=(\mathcal{U}^{\vec{x}},\mathcal{U}^{\vec{z}_0},\dots,\mathcal{U}^{\vec{z}_k},\mathcal{U}^{\vec{y}})$ and for every $j\leq k$ let $\mathbb{K}_k=(\mathcal{U}^{\vec{x}},\mathcal{U}^{\vec{z}_0},\dots,\mathcal{U}^{\vec{z}_k})$. We wil prove by induction that for every $i\leq n$, $\Gamma_{\mathbb{K}}(U_i^{\vec{x}})=U_i^{\vec{y}}$.\\\\
If $i=0$, notice that since $\{x_0,\dots, x_n\}\longrightarrow x_0$, then $\{U^{\vec{x}}_0,\dots U^{\vec{x}}_n\}\rightrightarrows U^{\vec{x}}_0$, which, by Claim 1, implies  that
$\{\Gamma_{\mathbb{K}_0}(U^{\vec{x}}_0),\dots \Gamma_{\mathbb{K}_0}(A^{\vec{x}}_n)\}\rightrightarrows 	\Gamma_{\mathbb{K}_0}(U^{\vec{x}}_0)$. Applying repeatedly Claim 1, and ussing the fact that for every $j\leq k$ if happens that $\Gamma_{(\mathcal{U}^{\vec{z}_{k}},\mathcal{U}^{\vec{z}_{k+1}})}\circ\Gamma_{\mathbb{K}_j}=\Gamma_{\mathbb{K}_{k+1}}$, we conclude that $$\{\Gamma_{\mathbb{K}}(U^{\vec{x}}_0),\dots \Gamma_{\mathbb{K}}(U^{\vec{x}}_n)\}\rightrightarrows 	\Gamma_{\mathbb{K}}(U^{\vec{x}}_0),$$
but $\Gamma_{\mathbb{K}}$ is a bijection  $\mathcal{U}^{\vec{x}}$ and $\mathcal{U}^{\vec{y}}$, so $\{\Gamma_{\mathbb{K}}(U^{\vec{x}}_0),\dots \Gamma_{\mathbb{K}}(U^{\vec{x}}_n\})=\{U^{\vec{y}}_0,\dots U^{\vec{y}}_n\}$. Thus, $\Gamma_{\mathbb{K}}(U^{\vec{x}}_0)=U^{\vec{y}}_0$.\\\\
Let $n>i>0$ and suppos that for every $j< i$ we had proven that $\Gamma_{\mathbb{K}}(U_j^{\vec{x}})=U_j^{\vec{y}}$.Proceeding in the same way as Case $i=0$, we can conclude that $$\{\Gamma_{\mathbb{K}}(U^{\vec{x}}_i),\dots \Gamma_{\mathbb{K}}(U^{\vec{x}}_n)\}\rightrightarrows 	\Gamma_{\mathbb{K}}(U^{\vec{x}}_i),$$
but $\Gamma_{\mathbb{K}}$ is a bijection from $\mathcal{U}^{\vec{x}}$ to $\mathcal{U}^{\vec{y}}$. Thus,by induction hypothesis  we have that $\{\Gamma_{\mathbb{K}}(U^{\vec{x}}_i),\dots \Gamma_{\mathbb{K}}(U^{\vec{x}}_n)\}=\{U^{\vec{y}}_i,\dots U^{\vec{y}}_n\}$. we conclude that $\Gamma_{\mathbb{K}}(U_i^{\vec{x}})=U_i^{\vec{y}}.$ With this, the induction ends.\\

To finish, notice that what we have just proved, implies that for every $\mathbb{K}$ and $\mathbb{L}$ $\mathfrak{D}$-chains from $\mathcal{U}^{\vec{x}}$ a $\mathcal{U}^{\vec{y}}$ we have that $\Gamma_{\mathbb{K}}=\Gamma_{\mathbb{L}}$. Thus, we have proved $a)\Rightarrow b)$.\\\\
\end{enumerate}

{$\bf{ \big(b)\Rightarrow a)\big)}$} For this, we only have to show that there exists a continuous selection over $\mathcal{P}(Reg_{n+1})$.\\\\
We define the relation $\sim\subseteq \mathcal{P}(Reg_{n+1})^2$ fiven \begin{center}$\vec{x}\sim \vec{y}$ if and only if there are  $\mathcal{U},\mathcal{V}\in \mathfrak{D}\; \mathfrak{D}$-chainables, such that  $\vec{x} \in \langle \mathcal{U}\rangle$  and  $\vec{y}\in \langle \mathcal{V}\rangle$.
 \end{center}	
 
It is easy to see that $\sim$ is an equivalence relation, which generates a partition $L(\sim)$ over $\mathcal{P}(Reg_{n+1})$. Given $C\in L(\sim) $ and $\vec{x}\in C$ there is $\mathcal{U}\in \mathfrak{D}$ such that $\vec{x}\in \langle \mathcal{U} \rangle$, and by definition $\sim$ we also have that $\langle\mathcal{U}\rangle\subseteq C$. By this observation, we can conclude that the elements of $L(\sim)$ are clopen over $\mathcal{P}(Reg_{n+1})$. Thus, we only need to show that for every $C\in L(\sim)$ there exists a continuous selection over $C$.\\

  Let $C\in L(\sim)$. For every $\vec{x}\in C$, fix $\mathcal{U}^{\vec{x}}\in \mathfrak{D}$ such that $\vec{x}\in \langle \mathcal{U}^{\vec{x}} \rangle$. Now fix $\vec{z}_0\in  C$, $U\in \mathcal{U}^{\vec{z}_0}$ and for every $\vec{x}\in C $ le'ts take  $\mathfrak{D}$-chain, namely, $\mathbb{K}_{\vec{x}}$ from $\mathcal{U}^{\vec{z}}_0$ to $\mathcal{U}^{\vec{x}}$. We define $f_C:C\longrightarrow X$ given by \begin{center} $f_C(\vec{x})$ is the unique element $\vec{x}$ in $\Gamma_{\mathbb{K}_{\vec{x}}}(U).$

 \end{center}

$f_C$ is well defined. Since $\mathfrak{D}$ is nice, for every $\vec{x}\in C$ and  $\vec{y}\in \langle \mathcal{U}^{\vec{x}}\rangle$ it happens that $f_C(\vec{y})\in \Gamma_{\mathbb{K}_{\vec{x}}}(U)$.  Thus, $f_C$ is continuous.

\end{proof}
\end{theorem}

\appendix





\bibliographystyle{elsarticle-num}

\end{document}